\magnification 1200\parindent 0pt\centerline
{\bf New Kazhdan groups \rm\footnote{$^1$}  
{both authors were partially supported by the KBN grant 2 P03A 023 14.
The first author was French government fellow under CIES bourse 231404E.}}

\vskip1cm
\centerline{Jan Dymara}
\medskip
\centerline{Instytut Matematyczny  Polskiej Akademii Nauk}
\centerline{Kopernika 18, 51--617 Wroc\l aw, Poland}
\centerline{\tt 
dymara@math.uni.wroc.pl}
%\centerline{http://www.math.uni.wroc.pl/\~{}dymara} 
\vskip 0.8cm
\centerline{Tadeusz Januszkiewicz}
\medskip
\centerline{Instytut Matematyczny, Uniwersytet Wroc\l awski}
\centerline{pl. Grunwaldzki 2/4, 50--384 Wroc\l aw, Poland}
\centerline{and}
\centerline{Instytut Matematyczny  Polskiej Akademii Nauk} 
\centerline{Kopernika 18, 51--617 Wroc\l aw, Poland}
\centerline{\tt tjan@math.uni.wroc.pl}
%\centerline{http://www.math.uni.wroc.pl/\~{}tjan}

\vskip 1cm
\centerline{\bf Introduction}

\bigskip
A locally compact, second countable group $G$ is called Kazhdan
if for any unitary representation of $G$, 
the first continuous cohomology group is trivial $H^1_{ct}(G, \rho)=0$.
There are several other equivalent definitions, the reader should consult 
[6], esp. 1.14 and 4.7.

For some time now Kazhdan groups have attracted attention.
One of the main challenges is to understand them geometrically.

Recently 
Pansu [7], \.Zuk [10], and Ballmann--\'Swi\c atkowski [1], went back to 
Garland's paper [5], improved it in several respects and produced 
among other things new examples of Kazhdan's groups. These examples, 
especially those in [1], are explicit and significantly different from 
classical ones. 

We also go back to Garland's paper, but instead of euclidean buildings,
we study hyperbolic ones.

\bigskip
An interesting class of
hyperbolic buildings with cocompact groups of automorphisms were constructed 
by Tits [9]. He associates with a ring $\Lambda$ and a generalised Cartan 
matrix $M$ a Kac--Moody group. These groups provide $BN$ pairs for buildings.
A special case of particular interest to us is that of $\Lambda$ a finite 
field and generalised Cartan matrix coming from hyperbolic reflection groups 
for which the fundamental domain is a simplex; there are 10 of them in 
dimension 2, two in dimension 3 and one in dimension 4.
Buildings associated to these data are locally finite and their 
automorphism groups are locally compact topological groups.
It turns out that they are Kazhdan' (and more).

\bigskip\bf Theorem 1. \sl 
Let $X_q$ be an $n$-dimensional building of thickness $(q+1)$, 
associated to a cocompact hyperbolic group with the fundamental domain 
a simplex.
Suppose $G$ is a closed in the compact open 
topology, unimodular subgroup of the simplicial 
automorphism group which
acts cocompactly on the building. 

Then 
for large $q$ and $1\leq k\leq n-1$
$$H^k_{ct}(G, \rho)=0,$$

that is the 
continous cohomology groups of $G$ with coefficients in any unitary 
representation vanish.
In particular $G$, considered as a topological group is a Kazhdan group.
\rm 

\bigskip
Several comments are in order:

1. 
The theorem holds for any hyperbolic building.
However at present Tits'  Kac--Moody buildings are the only examples where
we can verify the assumptions.

2.
For Tits' Kac--Moody building, 
simplicial automorphism groups which are \it uncountable, \rm are bigger than 
Kac--Moody groups (given by countably many generators and relations), 
and Tits' Kac--Moody groups are not discrete as subgroups of automorphism 
groups. In [8], B. R\'emy using twin buildings exhibits
Kac--Moody groups as discrete cofinite volume groups acting on
product of buildings, and also constructs discrete cofinite volume
groups acting on the building itself. His examples are not cocompact.

3. Unimodularity, brought in by the topology of the group, is an essential 
assumption. Kazhdan groups are unimodular.  On the other hand 
\'Swi\c atkowski pointed out to us nonunimodular groups acting cocompactly 
on classical euclidean buildings:
the upper triangular subgroup of $SL_n(Q_p)$ acting on its building.
In Tits' Kac--Moody examples, it is easy to establish that the group of
simplicial automorphisms is unimodular.

4. Three and four dimensional dimensional buildings provide first examples 
of Kazhdan groups of large dimension, not coming form locally symmetric 
spaces or euclidean buildings (they are also not products of lower 
dimensional examples, since they are hyperbolic).
Here "dimension" may be understood either as "continuous cohomological
dimension" or as "large scale dimension".
The argument here requires the computation of
$H^n_{ct}(G, {\rm St})$, where ${\rm St}$ is the
Steinberg representation of $G$ on the space of $l_2$-harmonic 
$n$-cochains on the building, computation which is not
essentially not different form the euclidean building case.
If we have discrete subgroups of automorphisms of these buildings,
they are necessarily Gromov hyperbolic, and the "dimension"
is one more than the dimension of their Gromov boundaries.

5. Bourdon [4] noticed that several two dimensional hyperbolic 
buildings admit cocompact actions of
discrete groups and thus one can use results of [1], 
[10] to show that some of them are Kazhdan. 
He does not use Tits' construction, but builds
his buildings as complexes of groups. Most of his buildings are not
Kac--Moody.

\bigskip
There are two ingredients we use in the proof. 

First is the Garland's method (which we take from [1], but actually [5]
implicitely contains almost all we need) 
for proving vanishing of cohomology groups of a 
simplicial complex.
Second is the use of continuous cohomology of topological groups,
in particular of Borel--Wallach result relating the cohomology of a 
complex on which the topological group acts with compact stabilizers, 
to its continuous cohomology.

The progress we obtain is that one does not have to worry about the 
existence of discrete subgroups. This is very handy, since bare existence 
of Tits examples is nontrivial, let alone their subtle properties.

\bigskip
In a future paper we construct more examples
of Kazhdan groups, all related to buildings.
\bigskip
{\sl Acknowledgements.} We are grateful to Jacek \'Swi\c atkowski 
for many useful discussions.

\bigskip\centerline{\bf \S 1 Generalities about automorphism groups.}\bigskip
We recall basic facts about simplicial automorphism groups
of simplicial complexes. They are all fairly standard.

\bigskip
Let $X$ be a countable locally finite simplicial complex, 
let $Aut(X)$ be the group of its simplicial automorphisms. 
The compact--open topology on $Aut(X)$ is defined 
using the basis of open neighbourhoods of the identity 
$U(K)=\{g: g\vert_K = id_K \}$, where $K$ runs over compact subsets in $X$. 

\bigskip
Let $G$ be a closed subgroup of $Aut(X)$, with induced topology.
Since $X$ is a countable complex, $Aut(X)$, thus $G$, has countable basis 
and hence it is metrizable by a left invariant metric.

\bigskip
\bf Proposition 1.1. \sl 
$G$ is locally compact. 
In fact stabilizers of compact subcomplexes in $X$ are compact and open. 
\rm 

\bigskip \bf Proposition 1.2. \sl $G$ is separable. \rm 

Thus, being metrizable and separable, $G$ is second countable
(has a countable basis).

\bigskip \bf Proposition 1.3. \sl $G$ is countable at infinity, 
or $\sigma$-compact: the sum of countably many compact subsets. \rm

\bigskip \bf Proposition 1.4. \sl Stabilizers of compact subcomplexes are 
either all finite or all uncountable. \rm

\bigskip \bf Proposition 1.5. \sl
$G$ is totally disconnected. \rm

\bigskip{\bf Unimodularity} 
of $G$ will play important role. Observe first that 
if the group $G$ is generated by compact subgroups then it is unimodular, 
since all generators go to identity under the modular homomorphism.

Suppose that a subgroup $G\subset Aut (X)$
generated by compact subgroups of $Aut (X)$ acts transitively on 
$n$-simplices. Then $Aut (X)$ is unimodular, since it is generated by 
$G$ and a stabilizer of a simplex.
A situation of interest to us where this happens is this:

\bigskip {\bf Lemma 1.6. }\sl 
Suppose $X$ is a connected locally finite simplicial complex, 
and suppose that links of simplices of codimension $\geq 2$ in $X$ 
are connected.
Suppose that stabilizers of $(n-1)$-simplices act transitively on their 
respective links.
Then $Aut(X)$ acts transitively on $X$ and is unimodular. \rm

The proof is clear from the above discussion. 
Observe that locally finite buildings coming from $BN$ pairs
satisfy the assumptions of the Lemma.

\bigskip\centerline {\bf \S 2 Borel--Wallach Lemma}\bigskip
Assume that $G$, a closed subgroup of the group of simplicial 
automorphisms of $X$, acts cocompactly. 
Sometimes one can identify $H^*_{ct}(G, \rho)$ with the 
cohomology of $X$ with coefficients in $\rho$.  Specifically:

Consider all alternating maps $\phi$
from ordered $k$-simplices in $X$ to $\cal H$,
satisfying  for all $g\in G$ and $\sigma \in X$

$$\phi (g\sigma)= \rho(g)\phi (\sigma)$$

Call the space of such maps $C^k(X,\rho)$. 
There is a natural differential making 
$C^*(X,\rho)$ into a complex

$$d\phi(\sigma)=\sum_{i=0}^k (-1)^i\sigma_i$$

where 
$\sigma=(v_0,\dots,v_k)$ and 
$\sigma_i=(v_0,\dots,v_{i-1}, v_{i+1},\dots,v_k)$

\bigskip
{\bf Lemma 2.1 } ([2], lemma X.1.12 page 297)
\sl
Let $(X, G)$ be an acyclic locally finite complex with a 
cocompact action of a group of its simplicial automorphisms.
Suppose $\rho$ is a representation of $G$ on a quasi-complete
(for example Hilbert) space. Then
$$H^i_{ct}(G, \rho)= H^i (C^*(X, \rho)).$$ \rm

\bigskip

The assumptions of this theorem are satisfied for locally finite 
buildings coming from $BN$ pairs.

\bigskip\centerline{\bf \S3 Vanishing theorem}\bigskip
Here we adapt Ballmann--\'Swi\c atkowski presentation of Garland's
method to greater generality of not necessarily discrete group actions.
To keep the exposition short we refer to their paper
for the notation.

\bigskip
{\bf Theorem 3.1.}
\sl Let $X$ be a locally finite simplicial complex,
and $G$ a cocompact unimodular group of its simplicial automorphisms.
Assume that for any simplex $\tau$ of $X$ the link $X_\tau$ is connected
and
$$\kappa_\tau>{k(n-k)\over k+1}$$ where $\kappa_\tau$
is the smallest positive eigenvalue of the Laplacian
$\Delta_\tau$ on $C^0(X_\tau, R)$.
\medskip
Then for $1\leq k\leq n-1$, 
$H^k(C^*(X,\rho))=0$ for any unitary representation $\rho$ of $G$.
\rm 

\bigskip
>From Lemma 2.1. we immediately get.
\medskip
{\bf Corollary 3.2.} 
\sl Under the assumptions of Theorem 3.1,
$G$ is a Kazhdan group, provided $X$ is acyclic.
\rm \medskip

Proof:
Theorem 3.1 corresponds to Theorem 2.5 of [1].
Their calculation goes through as it stands, except for two changes.

1. $\vert G_\sigma\vert$, there cardinality of the stabilizer
of $\sigma$, should now be understood as the Haar measure
of that stabilizer inside $G$.

2. Their Lemma 1.3 should be shifted form 
the discrete to locally compact setting.
Here is how this can be done.

Let $\Sigma (k)$
denote the set of ordered $k$-simplices in $X$,
let $\Sigma (k, G)$ be a set of representatves
of $G$  orbits on $\Sigma (k)$. Modified Lemma 1.3 of [1]
reads now as follows:

\medskip
{\bf Lemma 3.3.} \sl 
Let $X$, $G$ ba as in Theorem 3.1.
For $0\leq l< k\leq n$, let  $f=f(\tau,\sigma)$
be a $G$-invariant function on the set of pairs 
$(\tau,\sigma)$, where $\tau$ is an ordered $l$-simplex
and $\sigma$ is an ordered $k$-simplex
with $\tau\subset \sigma$, that is vertices of $\tau$
are vertices of $\sigma$. 
Then

$$
\sum_{\sigma \in \Sigma (k, G)} \
\sum_{\tau\in \Sigma (l) \atop \tau\subset\sigma} \
{f(\tau,\sigma)\over\vert G_\sigma\vert}
=
\sum_{\tau \in \Sigma (l, G)} \
\sum_{\sigma\in \Sigma (k)\atop \tau\subset \sigma} \
{f(\tau,\sigma)\over\vert G_\tau\vert}.
$$\rm 

Proof:
$$\eqalign{
\sum_{\sigma \in \Sigma (k, G)}
\sum_{\tau\in \Sigma (l) \atop
\tau\subset\sigma}
{f(\tau,\sigma)\over\vert G_\sigma\vert}
&=
\sum_
{\sigma\in\Sigma(k, G) \atop \tau\in\Sigma(l, G)}
\sum_{\gamma_i: 
\gamma_i\tau\subset\sigma \atop
\gamma_i\tau\neq\gamma_j\tau}
{f(\gamma_i\tau,\sigma)\over\vert G_\sigma\vert}
\cr&= 
\sum_{\sigma\in\Sigma(k, G) \atop \tau\in\Sigma(l, G)}
\int_{\bigcup_i\gamma_i G_\tau}
{f(\gamma\tau,\sigma)
\over
{\vert G_\tau\vert \vert G_\sigma\vert}}
d\gamma 
\cr&=
\sum_{\sigma\in\Sigma(k, G) \atop \tau\in\Sigma(l, G)}
\int_{\gamma:\gamma\tau\subset\sigma}
{f(\tau,\gamma^{-1}\sigma)\over
{\vert G_\tau\vert \vert G_\sigma\vert}}
d\gamma 
\cr&=
\sum_ {\sigma\in\Sigma(k, G) \atop \tau\in\Sigma(l, G)}
\int_{\gamma:\tau\subset\gamma\sigma}
{f(\tau,\gamma\sigma) \over
{\vert G_\tau\vert \vert G_\sigma\vert}}
d\gamma 
\cr&=
\sum_{\tau\in \Sigma (l,G)}
\sum_{\sigma \in \Sigma (k) \atop \tau\subset\sigma}
{f(\tau,\sigma)\over\vert G_\tau\vert} 
} $$

\bigskip\centerline{\bf \S4 Hyperbolic buildings}\bigskip

Here we rely on results of Tits [9].
We need the existence of buildings with cocompact proper group action,
with unimodular automorphism group, and arbitrarily large thickness.

Tits provides us with what we need as follows.
Take a Coxeter group, whose Dynkin diagram is 
a triangle, a square or a pentagon.
For triangles we allow labels $m,n,k$ on edges, such that
$m,n,k=2,3,4,6$, and ${1\over m}+{1\over n}+{1\over k}<1$.
For a square one of the edges is labelled 4, and remaining ones are 
labelled 3, or two opposite edges are labelled 4 and remaining ones are 
labelled 3 .
For pentagon one of the edges is labelled 4 and remaining ones are 
labelled 3.
These are all cocompact hyperbolic reflection group with 
the fundamental domain a simplex, and edges labelled $2,3,4,6$ [3].

For each such diagram and a finite field $F_q$, Tits constructs a Kac-Moody 
group, acting cocompactly (in fact transitively on simplices of maximal 
dimension)
on a hyperbolic building, with links of vertices being spherical
buildings of thickness $q+1$ corresponding to parabolic subgroups of the 
Coxeter system. Moreover the group is generated by 
elements stabilizing codimension 1 simplices.
Thus taking the closure of the Kac--Moody group
in the full automorphism group we obtain a cocompact unimodular
group acting on a hyperbolic building.

As far as we know, the existence of discrete cocompact subgroups 
in these groups has not been
established except for some two dimensional examples.

\bigskip
Now all we have to do to finish the proof of the theorem 1
is to check that the spectral condition holds
for the links. But this has been done (for large thickness)
already by Garland  [5, Sections 6--8] (see also remark at the 
end of the Section 3.1 of [1]).

\bigskip
It seems to us that Garland could have included these hyperbolic
examples in his original paper.

\vskip 1cm
\centerline{\bf Bibliography}
\bigskip

[1] Ballmann W. and \'Swi\c atkowski J.,
On $L^2$-cohomology and property (T) for automorphism groups of 
polyhedral cell complexes, 
GAFA  7(1997) 615--645

[2] Borel A. and Wallach N.,
Continuous cohomology, discrete subgroups, 
and representation of reductive groups,
Ann. of Math. Studies 94, Princeton University Press
and University of Tokyo Press 1980.

[3] Bourbaki N., 
"Groupes et algebres de Lie, chapitres IV-VI" Hermann 1968.

[4] Bourdon M., 
Sur les immeubles fuchsiennes et leur type de quasiisometrie
Preprint, Nancy, December 1997.

[5] Garland H.,
$p$-adic curvature and the cohomology of discrete subgroups
of $p$-adic groups,
Ann of Math. 97 (1973) 375--423.

[6] de la Harpe P.,  Valette A.,
La propri\'et\'e (T) de Kazhdan pour les groupes localement compacts,
Ast\'erisque 175, Soc. Math. France 1989.

[7]  Pansu P.,
Formule de Matsushima, de Garland, et propir\'et\'e (T) pour des groupes 
agissant sur des espaces symmetriques ou des immeubles,
Bull. Soc. Math. France 126(1998) pp. 107--139.

[8] Remy B.
Immeubles \`a courbure n\'egative et th\'eorie de Kac--Moody,
preprint (Nancy 1998).

[9] Tits J., 
Uniqueness and Presentation of Kac--Moody groups over fields,
J. of Algebra 105(1987) pp. 542--573.

[10] \.Zuk A., 
La propri\'et\'e (T) de Kazhdan pour les groupes agissant sur les poly\`edres, 
C. R. Acad. Sci. Paris 323, Serie I (1996), 453--458. 
\bye